\documentclass{article}

\usepackage{amsthm}
\usepackage{amsfonts}
\usepackage{amssymb}
\usepackage{amsgen}
\usepackage{amsmath}
\usepackage{amsopn}
\usepackage{verbatim}
\usepackage{xypic}
\usepackage{xspace}
\usepackage{multicol}
\usepackage{eepic}
\usepackage{upref}

\DeclareMathOperator{\im}{{Im}}

\DeclareMathOperator{\dch}{{dch}}
\DeclareMathOperator{\Lie}{{Lie}}

\DeclareMathOperator{\reg}{{reg}}
\DeclareMathOperator{\Gal}{{Gal}}

\theoremstyle{plain}
\newtheorem{thm}{Theorem}[section]
\newtheorem{lem}[thm]{Lemma}
\newtheorem{prop}[thm]{Proposition}
\newtheorem{cor}[thm]{Corollary}
\newtheorem{conj}[thm]{Conjecture}

\theoremstyle{definition}

\newtheorem{defn}[thm]{Definition}

\theoremstyle{remark}
\newtheorem{rem}[thm]{Remark}

\numberwithin{equation}{section}

 \font\cyr=wncyr10
 \newcommand{\nc}{\newcommand}
 \nc{\MPV}{{\mathcal{MPV}}}
 \nc{\calG}{{\mathcal G}}
 \nc{\wht}{{\widehat}}
 \nc{\bwg}{{\bigwedge}}
 \nc{\mmu}{{\boldsymbol{\mu}}}
 \nc{\mal}{{{\scriptstyle \maltese}}}
 \nc{\fA}{{\mathfrak A}}
 \nc{\HH}{{\mathfrak H}}
 \nc{\ra}{\rightarrow}
 \nc{\ors}{{\vec s\,}}
 \nc{\os}{{\overset}}
 \nc{\G}{{\mathbb G}}
 \nc{\Z}{{\mathbb Z}}
 \nc{\R}{{\mathbb R}}
 \nc{\N}{{\mathbb N}}
 \nc{\ZN}{{\mathbb Z_{\ge 0}}}
 \nc{\Q}{{\mathbb Q}}
 \nc{\C}{{\mathbb C}}
 \nc{\D}{{\mathcal D}}
 \nc{\calI}{{\mathcal I}}
 \renewcommand{\P}{{\mathbb P}}
 \nc{\caT}{{\mathcal T}}
 \nc{\tB}{{\tilde B}}
 \nc{\Li}{{\rm Li}}
 \nc{\suf}{{\ast\,}}
 \nc{\sufq}{{\ast_q\,}}
 \nc{\gam}{{\gamma}}
 \nc{\gG}{{\Gamma}}
 \nc{\om}{{\omega}}
 \nc{\ga}{{\alpha}}
 \nc{\gl}{{\lambda}}
 \nc{\gb}{{\beta}}
 \nc{\gd}{{\delta}}
 \nc{\gs}{{\sigma}}
 \nc{\gS}{{\Sigma}}
 \nc{\sif}{{\mathcal S}}
 \nc{\gt}{{\tau}}
 \nc{\Lra}{\Longrightarrow}
 \nc{\lra}{\longrightarrow}
 \nc{\fS}{{\mathfrak S}}
 \nc{\DD}{{\mathfrak D}}
 \nc{\Llra}{\Longleftrightarrow}
 \nc{\ol}{\overline}
 \nc{\lms}{\longmapsto}
 \nc{\cv}{{{\mathsf c}{\mathsf v}}}
 \nc{\zq}{{\zeta_q}}
 \nc\qup{{q\uparrow 1}}
 \nc{\us}{\underset}
 \nc{\tn}{{\tilde{n}}}
 \nc{\gD}{{\Delta}}
 \nc{\bk}{{\bf k}}
 \nc{\bi}{{\bf i}}
 \nc{\bfone}{{\bf 1}}
 \nc{\bfX}{{\bf X}}
 \nc{\bfY}{{\bf Y}}
 \nc{\QX}{{\Q\langle \bfX\rangle}}
 \nc{\QY}{{\Q\langle \bfY\rangle}}
 \nc{\CX}{{\C\langle \bfX\rangle}}
 \nc{\CY}{{\C\langle \bfY\rangle}}
 \nc{\QXX}{{\Q\langle\!\langle \bfX\rangle\!\rangle}}
 \nc{\QYY}{{\Q\langle\!\langle \bfY\rangle\!\rangle}}
 \nc{\CXX}{{\C\langle\!\langle \bfX\rangle\!\rangle}}
 \nc{\CYY}{{\C\langle\!\langle \bfY\rangle\!\rangle}}
 \nc{\bs}{{\bf s}}
 \nc{\caS}{{\mathfrak S}}
 \nc{\zz}{{\mathfrak z}}

\nc{\sha}{{\mbox{\cyr x}}}

\begin{document}

\title{Standard Relations
of Multiple Polylogarithm Values at Roots of Unity}
\author{Jianqiang Zhao}
\date{}
\maketitle
\begin{center}
{\large Department of Mathematics, Eckerd College, St. Petersburg,
FL 33711}
\end{center}
\vskip0.6cm

\noindent{\bf Abstract.}
Let $N$ be a positive integer. In this paper we shall study the
special values of multiple polylogarithms at $N$th roots of unity,
called multiple polylogarithm values (MPVs) of level $N$. These
objects are generalizations of multiple zeta values and
alternating Euler sums, which was studied by Euler, and more recently,
many mathematicians and theoretical physicists..
Our primary goal in this paper is to investigate the
relations among the MPVs of the same weight and level by using the
regularized double shuffle relations, regularized distribution relations,
lifted versions of such relations from lower weights,
and seeded relations which are produced by relations of weight one MPVs.
We call relations from the above four families \emph{standard}.
Let $d(w,N)$ be the $\Q$-dimension of $\Q$-span of all
MPVs of weight $w$ and level $N$. Then we obtain upper bound for
$d(w,N)$ by the standard relations which in general are no worse or
no better than the one given
by Deligne and Goncharov depending on whether $N$ is a prime-power or
not, respectively, except for 2- and 3-powers, in which case standard
relations seem to be often incomplete whereas Deligne shows that their
bound should be sharp by a variant of Grothedieck's period conjecture.
This suggests that in general there should
be other linear relations among MPVs besides the standard relations,
some of which are written down in this paper
explicitly with good numerical verification. We also provide a few
conjectures which are supported by our computational evidence.

\vskip0.6cm

\section{Introduction}
In recent years, there is a revival of interest in multi-valued
classical polylogarithms (polylogs) and their generalizations.
For any positive integers $s_1,\dots, s_\ell$,
Goncharov \cite{Gicm} defines the multiple polylogs of complex
variables as follows:
\begin{equation}\label{equ:polylog}
Li_{s_1,\dots, s_\ell}(x_1,\dots,x_\ell)=\sum_{k_1>\dots>k_\ell>0}
\frac{x_1^{k_1}\cdots x_\ell^{k_\ell}}{ k_1^{s_1}\cdots k_\ell^{s_\ell}}.
\end{equation}
Conventionally one calls $\ell$ the {\em depth} (or {\em length}) and
$s_1+\dots+s_\ell$ the {\em weight}. When the depth $\ell=1$ the
function is nothing but the classical polylog. When the weight is
also 1 we get the MacLaurin series of $-\log(1-x)$.
Another useful expression of the multiple polylogs is given
by the following iterated integral:
\begin{equation}\label{equ:iteratedForm}
Li_{s_1,\dots,s_\ell}(x_1,\dots,x_\ell)=(-1)^\ell
\int_{0}^1 \left(\frac{dt}{t}\right)^{\circ(s_1-1)}\circ\frac{d t}{t-a_1}
\circ\cdots \circ\left(\frac{dt}{t}\right)^{\circ(s_\ell-1)}\circ\frac{d t}{t-a_\ell}
\end{equation}
where $a_i=1/(x_1\dots x_i)$ for $1\le i\le \ell$. Here, we define the iterated
integrals recursively by $\int_a^b f(t)\circ w(t) =\int_a^b (\int_a^x f(t)) w(x)$
for any 1-form $w(t)$ and concatenation of 1-forms $f(t)$. We may think the
path lies in $\C$; however, it is more revealing to use iterated integrals
in $\C^\ell$ to find the analytic continuation of this function (see \cite{Zanamp}).

It is well-known that special values of polylogs have
significant applications in arithmetic such as Zagier's conjecture
\cite[p.622]{Zag}.  On the other hand, the multiple zeta
values (MZV) appear naturally in the study of the fundamental group
of $\P^1-\{0,1\infty\}$ which is closely related to the absolute Galois
group $\Gal(\overline{\Q}/\Q)$ according the Grothendieck \cite{D1}.
As pointed out by Goncharov, higher cyclotomy theory should study the
multiple polylogs at roots of unity, not only those of the classical ones.
Moreover, theoretical physicists have already found out
that such values appear naturally in the study of Feynmen diagrams
(\cite{Br1,Br2}).

Starting from early 1990s Hoffman \cite{H1,H2} has constructed
some quasi-shuffle (we will call ``stuffle'') algebras reflecting
the essential combinatorial properties of MZVs.
Recently he \cite{H3} extends this to
incorporate the multiple polylog values (MPVs) at roots of
unity, although his definition of $*$-product is different from
ours. Our approach here is a quantitative comparison between the
the results obtained by Racinet
\cite{Rac} who considers MPVs from the motivic viewpoint of Drinfeld
associators, and those by Deligne and Goncharov \cite{DG} who
study the motivic fundamental groups of $\P^1-(\{0,\infty\}\cup\mmu_N)$
by using the theory of mixed Tate motives over $S$-integers of number fields,
where $\mmu_N$ is the group of $N$th root of unity.

Fix an $N$th root of unity $\mu=\mu_N:=\exp(2\pi
\sqrt{-1}/N)$. The \emph{level} $N$ MPVs are defined by
\begin{equation}\label{equ:z}
L_N(s_1,\dots,s_\ell|i_1,\dots,i_\ell):=
Li_{s_1,\dots,s_\ell}(\mu^{i_1},\dots,\mu^{i_\ell}).
\end{equation}
We will always identify $(i_1,\dots,i_\ell)$ with $(i_1,\dots,i_\ell)
\pmod{N}.$ It is easy to see from \eqref{equ:polylog} that a MPV
converges if and only if  $(s_1,\mu^{i_1})\ne (1,1).$
Clearly, all level $N$ MPVs are
automatically of level $Nk$ for any positive integer $k$. For
example when $i_1=\cdots=i_\ell=0$ or $N=1$ we get the multiple
zeta values $\zeta(s_1,\dots,s_\ell)$. When $N=2$ we recover the
alternating Euler sums studied in \cite{BBB2,Zesum}. To save space,
if a substring
$S$ repeats $n$ times in the list then $\{S\}^n$ will be used. For
example, $L_N(\{2\}^2|\{0\}^2)=\zeta(2,2)=\pi^4/120$.

Standard conjectures in arithmetic geometry imply that $\Q$-linear relations
among MVPs can only exist between those of the same weight.
Let $\MPV(w,N)$ be the $\Q$-span of all the MPVs of weight $w$
and level $N$ whose dimension is denoted by $d(w,N)$.
In general, to determine $d(w,N)$ precisely is a very difficult
problem because any nontrivial lower bound would provide
some nontrivial irrational/transcendental results which is related to
a variant of Grothendieck's period conjecture (see \cite{Del}). For example, we
can easily show that $\MPV(2,4)=\langle \log^2 2, \pi^2, \pi \log 2\sqrt{-1},
(K-1)\sqrt{-1}\rangle,$ where $K=\sum_{n\ge 0} (-1)^n/(2n+1)^2$ is
the Catalan's constant. From Grothendieck's conjecture we know
$d(2,4)=4$ (see op. cit.) but we don't have a unconditional proof yet.
On the other hand, we may obtain upper bound of $d(w,N)$ by
finding as many linear relations in $\MPV(w,N)$ as possible.
As in the cases of MZVs and the alternating Euler sums
the regularized double shuffle relations (RDS) play
important roles in revealing the relations among MPVs. We shall
study this theory for MPVs in section \ref{sec:RDS}
by generalizing some results of \cite{IKZ} (also cf. \cite{Bigotte}).
It is commonly believed that in levels one and two all
linear relations among MPVs are consequences of RDS.

 From the point of view of Lyndon words and quasi-symmetric functions
Bigotte et al.\ \cite{Bigotte} have studied MPVs (they call them
\emph{colored MZVs}) primarily by using double shuffle relations.
However, when the level $N\ge 3$, these relations are not complete
in general, as we shall see in this paper.

If the level $N>3$ then by a theorem of Bass \cite{Bass} there are
many non-trivial linear relations (regarded as \emph{seeds})
in $\MPV(1,N)$ whose structure is clear to us.
Multiplied by MPVs of weight $w-1$
these relations can produce non-trivial linear relations
among MPVs of weight $w$ which we call the \emph{seeded
relations.} Similar to these relations we may produce new
relations by multiplying MPVs on RDS of lower weights. We call such
relations \emph{lifted relations}. We conjecture that when level
$N=3$ all linear relations among MPVs are consequences of the
RDS and the lifted RDS with $d(w,3)=2^w$.

Among MPVs we know that there are the so-called finite distribution
relations (FDT). Racinet \cite{Rac} considers further the regularization
of these relations by regarding MPVs as the coefficients of some
group-like element in a suitably defined pro-Lie-algebra of motivic origin.
Our computation shows that the regularized distribution relations (RDT) do
contribute to new relations not covered by RDS and FDT.
But they are not enough yet to produce all the lifted RDS.

\begin{defn}
We call a $\Q$-linear relation between MPVs \emph{standard}
if it can be produced by combinations of the following four families
of relations: regularized double shuffle relations (RDS),
regularized distribution relations (RDT), seeded
relations, and lifted relations from the above. Otherwise, it is
called a \emph{non-standard} relation.
\end{defn}

The main goal of this paper is to provide some numerical evidence
concerning the (in)completeness of the standard relations. Namely,
these relations in general are not enough
to cover all the $\Q$-linear relation between MPVs (see Remark~\ref{rem:Ta:bN}
and Remark~\ref{rem:wt2}); however, when the level is a prime $\le 47$ and
weight $w=2$ using a result of Goncharov we can show that the standard
relations are complete under the assumption of Grothendieck's period
conjecture (see \cite{Zocta}). We further find that when weight $w=2$
and $N=25$ or $N=49$, the standard relations are complete.
However, when $N$ is a 2-power or 3-power or
has at least two distinct prime factors, we know that the standard relations
are often incomplete by comparing our results with those of Deligne and
Goncharov \cite{DG}. Moreover, we don't know how to obtain the non-standard
relations except that when $N=4$,
we discover recently that octahedral symmetry of $\P^1-(\{0,\infty\}\cup \mmu_4)$
can produce some (presumably all) new relations not covered by the standard
ones (see op. cit.)

Most of the MPV identities in this paper are discovered with the
help of MAPLE using symbolic computations. We have verified almost
all relations by GiNaC \cite{GiNac} with an error bound $<10^{-90}$.

This work was started while I was
visiting Chern Institute of Mathematics at Nankai University and
the Morningside Center of Mathematics at
Beijing, China in the summer of 2007. I would like to thank both
institutions and my hosts Chengming Bai and Fei Xu for their
hospitality and the ideal working environment.
I also want to thank Jens Vollinga for answering some of my questions
regarding the numerical computation of the multiple polylog values.
The paper was revised later while I was visiting the
Institute for Advanced Study and thanks are due to
Prof. Deligne for his patient explanation
of \cite{DG} and many insightful remarks on the paper.
This work was partially supported by a faculty
development fund from Eckerd College.

\section{The double shuffle relations and the algebra $\fA$}\label{sec:RDSoverC}
It is Kontsevich \cite{K1} who first noticed that MZVs can be
represented by iterated integrals (cf.~\cite{Rac}). We now extend this to MPVs. Set
$$a=\frac{dt}{t},\qquad b_i=\frac{\mu^i dt}{1-\mu^i t}\quad
\text{ for } i=0,1,\dots,N-1.$$
For every positive integer $n$ define
$$y_{n,i}:=a^{n-1} b_i.$$
Then it is straight-forward to verify using \eqref{equ:iteratedForm}
that if $(s_1,\mu^{i_1})\ne (1,1)$ then (cf.~\cite[(2.5)]{Rac})
\begin{equation}\label{equ:mzv}
L_N(s_1,\dots,s_n|i_1,i_2,\dots,i_n)=\int_0^1
y_{s_1,i_1}y_{s_2,i_1+i_2}\cdots
 y_{s_n,i_1+i_2+\dots+i_n}.
\end{equation}
We now define an algebra of words as follows:
\begin{defn} Set $A_0=\{\bfone\}$ to be the set of the empty
word. Define $\fA=\Q\langle A\rangle$ to be the graded
noncommutative polynomial $\Q$-algebra generated by letters $a$
and $b_i$ for $i\equiv 0,\dots,N-1\pmod{N}$, where $A$ is a
locally finite set of generators whose degree $n$ part $A_n$
consists of words (i.e., a monomial in the letters) of depth $n$.
Let $\fA^0$ be the subalgebra of $\fA$ generated by words not
beginning with $b_0$ and not ending with $a$. The words in $\fA^0$
are called \emph{admissible words.}
\end{defn}

Observe that every MPV can be expressed uniquely as an iterated
integral over the closed interval $[0,1]$ of an admissible word
$w$ in $\fA^0$. Then we denote this MPV by
\begin{equation}\label{equ:Z}
Z(w):=\int_0^1 w.
\end{equation}
Therefore we have (cf.~\cite[(2.5) and (2.6)]{Rac})
\begin{align}
 \label{equ:1-1LZ}
L_N(s_1,\dots,s_n|i_1,i_2,\dots,i_n)=&Z(
y_{s_1,i_1}y_{s_2,i_1+i_2}\cdots
 y_{s_n,i_1+i_2+\dots+i_n}),\\
 Z(y_{s_1,i_1}y_{s_2,i_2}\cdots
 y_{s_n,i_n})=&L_N(s_1,\dots,s_n|i_1,i_2-i_1,\dots,i_n-i_{n-1}).
 \label{equ:1-1ZL}
\end{align}
For example $L_3(1,2,2|1,0,2) =Z(y_{1,1}y_{2,1}y_{2,0}).$
On the other hand, during 1960s Chen developed a theory of
iterated integral which can be applied in our situation.
\begin{lem}\label{chen's}
\emph{(\cite[(1.5.1)]{Chen})}
Let $\om_i$ $(i\ge 1)$ be $\C$-valued 1-forms on a manifold $M$. For
every path $p$,
 $$ \int_p \om_1\cdots \om_r\int_p \om_{r+1}\cdots \om_{r+s}=
 \int_p (\om_1\cdots \om_r) \sha (\om_{r+1}\cdots \om_{r+s})$$
where $\sha$ is the shuffle product defined by
 $$(\om_1\cdots \om_r) \sha (\om_{r+1}\cdots
 \om_{r+s}):=\sum_{\substack{\gs\in
 S_{r+s},\gs^{-1}(1)<\cdots<\gs^{-1}(r)\\
 \gs^{-1}(r+1)<\cdots< \gs^{-1}(r+s)}}
   \om_{\gs(1)}\cdots \om_{\gs(r+s)}.$$
\end{lem}
For example, we have
  \begin{align*}
\ &L_N(1|1)L_N(2,3|1,2)=Z(y_{1,1})Z(y_{2,1}y_{3,3})=Z(b_1\sha
(ab_1a^2b_3))\\
 =&Z(b_1ab_1a^2b_3+2ab_1^2a^2b_3+(ab_1)^2ab_3+ab_1a^2b_1b_3+ab_1a^2b_3b_1)\\
 =&Z(y_{1,1}y_{2,1}y_{3,3}+2y_{2,1}y_{1,1}y_{3,3}+y_{2,1}^2y_{2,3}
 +y_{2,1}y_{3,1}y_{1,3}+y_{2,1}y_{3,3}y_{1,1})\\
 =&L_N(1,2,3|1,0,2)+2L_N(2,1,3|1,0,2)+L_N(2,2,2|1,0,2)\\
 \ &\hskip3cm +L_N(2,3,1|1,0,2)+L_N(2,3,1|1,2,N-2).
\end{align*}

Let $\fA_\sha$ be the algebra of $\fA$ together with the
multiplication defined by shuffle product $\sha$. Denote the
subalgebra $\fA^0$ by $\fA_\sha^0$ when we consider the shuffle
product. Then we can easily prove
\begin{prop} \label{shahomo} The map $Z:  \fA_\sha^0\lra \C$
is an algebra homomorphism.
\end{prop}

On the other hand, it is well known that MPVs also satisfy the
series stuffle relations. For example
\begin{equation*}
   L_N(2|5)L_N(3|4)=L_N(2,3|5,4)+L_N(3,2|4,5)+L_N(5|9).
\end{equation*}
because
$$\sum_{j>0}\sum_{k>0}=\sum_{j>k>0}+\sum_{k>j>0}+\sum_{j=k>0}.$$
To study such relations in general we need the following
definition.
\begin{defn} Denote by $\fA^1$ the subalgebra of
$\fA$ which is generated by words $y_{s,i}$ with $s\in \Z_{>0}$
and $i\equiv 0,\dots,N-1\pmod{N}$. Equivalently, $\fA^1$ is the
subalgebra of $\fA$ generated by words not ending with $a$. For
any word $w=y_{s_1,i_1}y_{s_2,i_2}\cdots
 y_{s_n,i_n}\in \fA^1$ and positive integer $j$ we define the exponent shifting
operator $\tau_j$ by
 $$\tau_j(w)=y_{s_1,j+i_1}y_{s_2,j+i_2}\cdots y_{s_n,j+i_n}.$$
For convenience, on the empty word we have the convention that
$\tau_j(\bfone)=\bfone.$ We then define a new multiplication $*$
on $\fA^1$ by requiring that $*$ distribute over addition, that
$\bfone*w=w*\bfone=w$ for any word $w$, and that, for any words
$\om_1,\om_2$,
\begin{multline}  \label{equ:defnstuffle}
 y_{s,j}\om_1*y_{t,k}\om_2 = y_{s,j}\Big(\tau_j\big(\tau_{-j}(\om_1)*y_{t,k}\om_2\big)\Big)
 + y_{t,k}\Big(\tau_k\big( y_{s,j}\om_1*\tau_{-k}(\om_2)\big)\Big) \\
 +y_{s+t,j+k}\Big(\tau_{j+k}\big(\tau_{-j}(\om_1)*\tau_{-k}(\om_2)\big)\Big).
\end{multline}
We call this multiplication the \emph{stuffle product}.
\end{defn}

\begin{rem} Our $\fA$, $\fA^0$ and $\fA^1$ are related to
$\QX$, $\QX_\cv$ and $\QY$ of \cite{Rac}, respectively. See section \ref{sec:Rac}.
\end{rem}

If we denote by $\fA_*^1$ the algebra $(\fA^1,*)$ then it is not
hard to show that (cf.~\cite[Thm.~2.1]{H2})
\begin{thm}
The polynomial algebra $\fA_*^1$ is a commutative graded
$\Q$-algebra.
\end{thm}

Now we can define the subalgebra $\fA_*^0$ similar to $\fA_\sha^0$
by replacing the shuffle product by stuffle product. Then by
induction on the lengths and using the series definition we can
quickly check that for any $\om_1,\om_2\in \fA_*^0$
$$Z(\om_1)Z(\om_2)=Z(\om_1\ast \om_2).$$
This implies that
\begin{prop} \label{*homo}
The map $Z:  \fA_*^0  \lra \C$ is an algebra homomorphism.
\end{prop}

For $\om_1,\om_2\in \fA^0$ we will say that
 $$Z(\om_1\sha \om_2-\om_1*\om_2)=0$$
is a finite double shuffle (FDS) relation. It is known that even
in level one these relations are not enough to provide all the
relations among MZVs. However, it is believed that one can remedy this by
considering RDS produced by the
following mechanism. This was explained in detail in \cite{IKZ}
when Ihara, Kaneko and Zagier considered MZVs where they call these
extended double shuffle relations.

Combining Propositions \ref{*homo} and \ref{shahomo} we can prove
easily (cf.~\cite[Prop.~1]{IKZ}):
\begin{prop} \label{prop:eDS}
We have two algebra homomorphisms:
$$Z^*: (\fA_*^1,*)\lra \C[T],\quad \text{and}\quad Z^\sha: (\fA_\sha^1,\sha)\lra \C[T]$$
which are uniquely determined by the properties that they both
extend the evaluation map $Z:\fA^0\lra \C$ by sending
$b_0=y_{1,0}$ to $T$.
\end{prop}

In order to establish the crucial relation between $Z^*$ and
$Z^\sha$ we can adopt the machinery in \cite{IKZ}. For any
$(\bs|\bi)=(s_1,\dots,s_n|i_1,\dots,i_n)$ where $i_j$'s are
integers and $s_j$'s are positive integers, let the image of the
corresponding words in $\fA^1$ under $Z^*$ and $Z^\sha$ be denoted
by $Z_{(\bs|\bi)}^*(T)$ and $Z_{(\bs|\bi)}^\sha(T)$ respectively.
For example,
\begin{align*}
 TL_N(2|3)=&Z_{(1|0)}^*(T)Z_{(2|3)}^*(T)=Z^*(y_{1,0}*y_{2,3})\\
 =&Z_{(1,2|0,3) }^*(T)+Z_{(2,1|3,3) }^*(T)+Z_{(3|3)}^*(T),
\end{align*}
while
\begin{align*}
 TL_N(2|3)=&Z_{(1|0)}^\sha(T)Z_{(2|3)}^\sha(T)=Z^\sha(y_{1,0}\sha
 y_{2,3})=Z^\sha(b_0\sha ab_3)\\
 =&Z_{(1,2|0,3)}^\sha(T)+Z_{(2,1|0,3)}^\sha(T) +Z_{(2,1|3,0)}^\sha(T)  .
\end{align*}
Hence we find the following RDS by the next Theorem:
 $$L_N(2,1|3,0) +L_N(3|3)=L_N(2,1|3,N-3)+L_N(2,1|0,3).$$
\begin{thm}\label{thm:RDSoverC}
Define a $\C$-linear map $\rho:\C[T]\to \C[T]$ by
 $$\rho(e^{Tu})=\exp\left(\sum_{n=2}^\infty\frac{(-1)^n}{n}\zeta(n)u^n\right)e^{Tu},\qquad |u|<1.$$
Then for any index set $(\bs|\bi)$ we have
 $$Z_{(\bs|\bi)}^\sha(T)= \rho(Z_{(\bs|\bi)}^*(T)).$$
\end{thm}
This is a the generalization of \cite[Thm.~1]{IKZ} to the
higher level MPV cases. The proof is essentially the same. One may
compare Cor.~2.24 in \cite{Rac}. The above steps  can be easily
transformed to computer codes which are used in our MAPLE programs.

\section{Finite and regularized double shuffle relations (FDS \& RDS)}\label{sec:RDS}
It is generally believed that all the linear relations between MZVs
can be derived from RDS.
Although the naive generalization of this to arbitrary levels is wrong the
idea in \cite{IKZ} to formalize this via some universal objects is still very useful.
We want to generalize this idea to MPVs in this section.

Keep the same notation as in the preceding sections. Let $R$ be a commutative
$\Q$-algebra with 1 and $Z_R:\fA^0\lra R$ such
that the ``finite double shuffle'' (FDS) property holds:
  $$Z_R(\om_1\sha \om_2)=Z_R(\om_1* \om_2)=Z_R(\om_1)Z_R(\om_2).$$
We then extend $Z_R$ to $Z_R^\sha$ and $Z_R^*$ as before. Define
an $R$-module $R$-linear automorphism $\rho_R$ of $R[T]$ by
$$\rho_R(e^{Tu})=A_R(u)e^{Tu}$$
where
$$A_R(u)=\exp\left(\sum_{n=2}^\infty \frac{(-1)^n}{n}
 Z_R(a^{n-1} b_0)u^n\right) \in R[\![u]\!].$$
Similar to the situation for MZVs, we may define the
$\fA^0$-algebra isomorphisms
$$\reg_\sha^T:\fA_\sha^1 =\fA_\sha^0[b_0]\lra\fA_\sha^0[T],\qquad
 \reg_*^T:\fA_*^1 =\fA_*^0[b_0]\lra\fA_*^0[T],$$
which send $b_0$ to $T$. Composing these with the evaluation map
$T=0$ we get the maps $\reg_\sha$ and $\reg_*$.
\begin{thm} \label{thm:RDS}
Let $(R,Z_R)$ be as above with the FDS property. Then the
following are equivalent:
\begin{itemize}
    \item[\emph{(i)}] $(Z_R^\sha-\rho_R\circ Z_R^*)(w)=0$ for all $w\in \fA^1$.
    \item[\emph{(ii)}]  $(Z_R^\sha-\rho_R\circ Z_R^*)(w)|_{T=0}=0$ for all $w\in \fA^1$.
   \item[\emph{(iii)}]  $ Z_R^\sha(\om_1\sha \om_0-\om_1*\om_0)=0$ for all $\om_1\in
   \fA^1$ and all $\om_0\in \fA^0$.
   \item[\emph{(iii$'$)}]  $ Z_R^*(\om_1\sha \om_0-\om_1*\om_0)=0$ for all $\om_1\in
   \fA^1$ and all $\om_0\in \fA^0$.
   \item[\emph{(iv)}]  $ Z_R(\reg_\sha(\om_1\sha \om_0-\om_1*\om_0))=0$ for all $\om_1\in
   \fA^1$ and all $\om_0\in \fA^0$.
   \item[\emph{(iv$'$)}]  $ Z_R(\reg_*(\om_1\sha \om_0-\om_1*\om_0))=0$ for all $\om_1\in
   \fA^1$ and all $\om_0\in \fA^0$.
   \item[\emph{(v)}]  $ Z_R(\reg_\sha(b_0^m*w))=0$ for all $m\ge 1$ and all $w\in \fA^0$.
   \item[\emph{(v$'$)}]  $ Z_R(\reg_*(b_0^m\sha w-b_0^m*w))=0$ for all $m\ge 1$ and all $w\in \fA^0$.
\end{itemize}
If $Z_R$ satisfies any one of these then we say that $Z_R$ has the
\emph{regularized double shuffle} (RDS) property.
\end{thm}
Notice that RDS automatically implys FDS.
The proof of the theorem is almost the same as that of
\cite[Thm.~2]{IKZ} but for completeness we give the most
important details in the following because there is some subtle
difference for MPVs of arbitrary level.

Denote by $\caS$ the set of the $y_{s,j}$ ($s\in \Z_{>0},
j=0,\dots, N-1$). For convenience we write $\tau_z=\tau_j$ if
$z=y_{s,j}\in\caS$. If $w=y_{s_1,i_1}\dots y_{s_n,i_n}\in \fA^1$
then we put $\tau_w=\tau_{i_1+\dots+i_n}$ and
$\tau_{-w}=\tau_{-i_1-\dots-i_n}$. Then (cf. \cite[Prop.~2]{IKZ})
\begin{prop} \label{prop:gdz} We have
\begin{enumerate}
    \item [(i)] For $z\in \caS$ the map $\gd_z:\fA^1\to \fA^1$ defined by
$$\gd_z(w):=z*w-z\tau_z(w)$$
is a ``twisted derivation'' in the sense that
$$\gd_z(ww')=\gd_z(w)\tau_z(w')+w\tau_w\Big(\gd_z\big(\tau_{-w}(w')\big)\Big).$$
Moreover, all these twisted derivations commute.
    \item [(ii)] The above twisted derivations extend to a twisted
    derivation on all of $\fA$ after setting $\tau_a=\text{id}$, with values on the letters $a,
    b_j$ given by
    $$\gd_z(a)=0, \quad \gd_z(b_j)=(a+b_j)\tau_j(z)\quad (z\in \caS, j=0,\dots,N-1).$$
    In particular, $\gd_z$ preserves $\fA^0$.
\end{enumerate}
\end{prop}
\begin{proof}  Easy computation by Definition
\eqref{equ:defnstuffle}.
\end{proof}
\begin{cor} \label{cor:gdz}
Denote by $\zz$ the $\Q$-linear span of the $y_{s,0}$ ($s\in
\Z_{>0}$). Then for $z\in \zz$ the map $\gd_z$ is a derivation on
$\fA^1$ which preserves $\fA^0$. Moreover, $\gd_z$ can be extended
to a derivation on $\fA$ by Prop.~\ref{prop:gdz}(ii).
\end{cor}
\begin{proof} Define $\gd'$ as in Prop.~\ref{prop:gdz}(ii). For any
$z=y_{s,0}\in\zz$ and $y_{t,i}$ a generator of $\fA^1$ we have
$$\gd'_z(y_{t,i})=a^{t-1}\gd'_z(b_i)=a^{t-1}(a+b_i)y_{s,i}
 =y_{s,0}*y_{t,i}-y_{s,0}y_{t,i}=\gd_z(y_{t,i}).$$
Note that for $z\in \zz$ and $w, w'\in\fA$ we have
$\tau_w\Big(\gd_z\big(\tau_{-w}(w')\big)\Big)=\gd_z(w')$. So
indeed $\gd_z$ can be extended to a derivation on $\fA$. It's
obvious that $\gd_z$ fixes both $\fA^0$ and $\fA^1$.
\end{proof}
We can define another operation on $\caS$ by $y_{s,i}\circ y_{t,j}=y_{s+t,i+j}.$
We can then restrict this to $\{y_{s,0}: s\in\Z_{>0}\}$ then
extend linearly to $\zz$. The following result is then
straight-forward.
\begin{prop} The vector space $\zz$ becomes a commutative and
associative algebra with respect to the multiplication $\circ$
defined by
$$z\circ z'=z*z'-zz'-z'z.$$
\end{prop}
The following proposition is one of the keys to the proof of
Theorem~\ref{thm:RDS} (cf. \cite[Prop.~4]{IKZ}).
\begin{prop}\label{circ}
Let $u$ be a formal variable. For $z\in \caS$ we have
$$\exp(zu\tau_z) (\bfone) =(2-\exp_{\circ}(z u\tau_z) )^{-1}  (\bfone).$$
(The inverse on the right is with respect to the concatenation
product.)
\end{prop}
\begin{proof} Define power series
$$f(u \tau_z)=\exp_\circ(z u\tau_z )-1=z u\tau_z +z\circ z \frac{u^2\tau_z^2}{2}+\cdots  $$
Then taking derivative with respect to $u$ we get
$f'(u\tau_z)=z\circ(1+f(u \tau_z))\tau_z.$
Now for $z,\om_i\in\caS$ we have by Prop.~\ref{prop:gdz} and Prop.~\ref{circ}
$$z*(\om_1\om_2\cdots \om_n)=\sum_{i=0}^n \om_1\cdots
\om_i\tau_{\om_i}(z)\tau_z(\om_{i+1}\cdots \om_n)
+\sum_{i=1}^n \om_1\cdots
 \om_{i-1}( z\circ \om_i)\tau_z(\om_{i+1}\cdots \om_n).$$
This yields
  \begin{align*}
\ & z*\left(z^{\circ n_1}\frac{(u \tau_z)^{n_1}}{n_1!}\cdots
 z^{\circ n_d}\frac{(u \tau_z)^{n_d}}{n_d!}\right)\\
=& \sum_{i=0}^d \left(z^{\circ n_1}\frac{(u
\tau_z)^{n_1}}{n_1!}\cdots
 z^{\circ n_i}\frac{(u \tau_z)^{n_i}}{n_i!}\right)z\tau_z \left(z^{\circ n_{i+1}}\frac{(u
\tau_z)^{n_{i+1}}}{n_{i+1}!}\cdots
 z^{\circ n_d}\frac{(u \tau_z)^{n_d}}{n_d!}\right)\\
+&\sum_{i=1}^d
 \left(z^{\circ n_1}\frac{(u \tau_z)^{n_1}}{n_1!}\cdots
 z^{\circ (n_i+1)}\frac{(u \tau_z)^{n_i}}{n_i!}\right)\tau_z\left(z^{\circ n_{i+1}}\frac{(u
\tau_z)^{n_{i+1}}}{n_{i+1}!}\cdots
 z^{\circ n_d}\frac{(u \tau_z)^{n_d}}{n_d!}\right).
\end{align*}
 Hence
 $$z *\big(1-f(u \tau_z)\big)^{-1}\tau_z(\bfone)=\frac{d}{du}\Big(\big(1-f(u
 \tau_z)\big)^{-1} \Big) (\bfone).$$
 This implies that
$$\exp_*(zu \tau_z ) (\bfone)=\big(1-f(u\tau_z)\big)^{-1} (\bfone)$$
as desired.
\end{proof}

\begin{cor}\label{cor:key} For all $z\in \caS$ we have
$$\exp_*(\log_\circ (1+z\tau_z)) (\bfone)=\big(1-z\tau_z\big)^{-1} (\bfone).$$
\end{cor}
If $z\in \zz$ then $\tau_z=\text{id}$ and therefore we have
\begin{cor}\label{cor:circ} For $z\in \zz$ we have
 \begin{equation*}
 \exp(zu)=(2-\exp_{\circ}(z u) )^{-1},\quad \ \exp_*(\log_\circ
(1+z))=(1-z)^{-1}.
\end{equation*}
\end{cor}

Let's consider a non-trivial example of Cor.~\ref{cor:key}. Let
$N=2$ and $z=y_{k,1}$ then we have
$$\exp\left(\sum_{n=1}^\infty (-1)^{n-1}
\zeta\big(nk;(-1)^k\big)\frac{u^n}{n}\right)=1+\sum_{n=1}^\infty
 \zeta\big(\{k\}^n;\{-1\}^n\big) u^n$$
where $\zeta(s_1,\dots,s_d;(-1)^{\gs_1},\dots,(-1)^{\gs_d})=L_2(
s_1,\dots,s_d|\gs_1,\dots,\gs_d)$ are the alternating Euler sums.
For example, by comparing the coefficients of $u^2$ and $u^3$ we get
$2\zeta(\ol{k},\ol{k})=\zeta(\ol{k})^2-\zeta(2k),$ and
$6\zeta(\ol{k},\ol{k},\ol{k})=\zeta(\ol{k})^3-3\zeta(\ol{k})\zeta(2k)+2\zeta(\ol{3k}).$
Here $\ol{s_j}$ means that the corresponding $\gs_j$ is odd.

The following two propositions are generalizations of Prop.~5-6 of
\cite{IKZ} respectively whose computational proofs are mostly
omitted since nothing new happens in our situation.

\begin{prop} \label{prop:zz'}
   For $z,z'\in \zz$ and $w\in \fA^1$ we have
\begin{align}
 \label{equ:gdzz'}
    \exp(\gd_z)(z')=&\big(\exp_\circ(z)\circ z'\big)\exp_*(z),\\
   \exp(\gd_z)(w)=&(\exp_*(z))^{-1}\big(\exp_*(z)* w\big).
\end{align}
\end{prop}

\begin{prop} \label{prop:Phi}
   For $z\in \zz$ define $\Phi_z:\fA^1\to \fA^1$ by
\begin{equation}\label{equ:Phi}
    \Phi_z(w):=(1-z)\left(\big(1-z\big)^{-1}*w\right)\qquad (w\in
    \fA^1).
\end{equation}
Then  $\Phi_z$ is an automorphism of $\fA^1$ and we have
\begin{equation}\label{equ:Phiexp}
 \Phi_z(w)=\exp(\gd_t)(w), \qquad\text{where }t=\log_\circ(1+z)\in
 \zz.
\end{equation}
All the $\Phi_z$ commute. Moreover, after restricting the
derivation $\gd_t$ to $\fA^0$ we can regard $\Phi_z$ as an
automorphism of $\fA^0$, If we extend the derivation $\gd_t$ to
the whole $\fA$ as in Cor.~\ref{cor:gdz} then we can regard
$\Phi_z$ as an automorphism of $\fA$.
\end{prop}
\begin{proof} The key point is that $\gd_t$ sends
$\fA^0$ to $\fA^0$ as a derivation by Cor.~\ref{cor:gdz}. Hence
$\exp(\gd_t)$ is an automorphism on $\fA^1$ as well as on $\fA^0$.
\end{proof}

The next three results are generalizations of Prop.~7, its
corollary, and Prop.~8 of \cite{IKZ}, respectively. The proofs
there can be easily adapted into our situation because the
$\sha$-product is essentially the same (note that the the only
essentially new phenomenon in the higher level MPV cases is that
there are exponent shiftings on the roots of unity in our stuffle
product.)
\begin{prop} \label{prop:d}
   Define the map $d:\fA\to \fA$ by $d(w)=b_0\sha w-b_0 w$. Then
   $d$ is a derivation and by setting $u$ as a formal parameter we have
\begin{equation*}
    \exp(du)(w)=(1-b_0u)\left(\big(1-b_0u\big)^{-1}\sha w\right)\qquad (w\in
    \fA^1).
\end{equation*}
On the generators we have
\begin{equation}\label{equ:dgens}
 \exp(du)(a)=a\big(1-b_0u\big)^{-1},\qquad\exp(du)(b_j)=b_j\big(1-b_0u\big)^{-1},\
 j=0,\dots,N-1.
\end{equation}
\end{prop}
\begin{rem}
In fact, we can replace the whole $\fA$ by $\fA^1$ in the first
part of Prop.~\ref{prop:d}. We can do the same in the next
corollary. However, in the proof of Theorem~\ref{thm:RDS} we only
need this weaker version.
\end{rem}

\begin{cor} \label{cor:Delta} Let $u$ be a formal parameter.
Let $\gD_u=\exp(-du)\circ \Phi_{b_0u}\in {\rm Aut}(\fA)[\![u]\!]$
(here $\circ$ means the composition). Then
 \begin{equation*}
 (1-b_0u)^{-1}*w=(1-b_0u)^{-1} \sha \gD_u(w), \qquad \forall
w\in\fA^1.
\end{equation*}
In particular, for $w\in \fA^0$ by taking $\reg$ on both sides of
the above equation we get
\begin{equation*}
    \reg_\sha\Big((1-b_0u)^{-1}*w\Big)=  \gD_u(w).
\end{equation*}
\end{cor}

\begin{prop} \label{prop:reg}
 For $\om_0=a\om_0' \in\fA^0$ we have
 \begin{equation*}
   \reg_\sha^T\left(\big(1-b_0u\big)^{-1} \om_0\right)=
   \exp(-du)(\om_0)e^{Tu}=a\left(\big(1+b_0u\big)^{-1}\sha \om_0'\right)e^{Tu}.
\end{equation*}
\end{prop}

\begin{rem}  Theorem~\ref{thm:RDS} now follows easily from
a detailed computation as in \cite{IKZ}. As a matter of fact, the
same argument shows that \cite[Prop.~10]{IKZ} and its Cor.~are
both valid in our general setup if we replace $\HH^0$ there by
$\fA^0$.
\end{rem}

\section{Seeded (or weight one) relations} \label{sec:seed}
When $N\ge 4$ there exist linear relations among MPVs of weight one by a
theorem of Bass \cite{Bass}. These relations are important because by
multiplying any MPV of weight $w-1$ by such a relation we can
get a relation between MPVs of weight $w$ which we call a
\emph{seeded relation.} This is one of the key ideas in finding
the formula in \cite[5.25]{DG} concerning $d(w,N)$.

First, we know there are $N-1$ MPVs of weight 1 and level $N$:
 $$L_N(1|j)=-\log(1-\mu^j), \qquad 0<j<N ,$$
where $\mu=\mu_N=\exp(2\pi\sqrt{-1}/N)$ as before. Here we have
taken $\C-(-\infty,0]$ as the principle domain of the logarithm.
Further, it follows from
the motivic theory of classical polylogs developed by Deligne
and Beilinson and the Borel's theorem (see \cite[Thm.~2.1]{G1}) that
the $\Q$-dimension of $\MPV(1,N)$ is
$$d(1,N)=\dim K_1(\Z[\mu_N][1/N])\otimes \Q+1=\varphi(N)/2+\nu(N),$$
where $\varphi$ is the Euler's totient function and $\nu(N)$
is the number of distinct prime factors of $N$.
Hence there are many linear relations
among $L_N(1|j)$. For instance, if $j<N/2$ then we have the
symmetric relation
$$-\log(1-\mu^j)=-\log(1-\mu^{N-j})-\log(-\mu^j)=
 -\log(1-\mu^{N-j})+\frac{N-2j}{N}\pi\sqrt{-1}.$$
Thus for all $1<j<N/2$
\begin{equation}\label{equ:symrel}
(N-2)(L_N(1|j)-L_N(1|N-j))=(N-2j)(L_N(1|1)-L_N(1|N-1)).
\end{equation}
Further, from \cite[(B)]{Bass} for any divisor $d$ of $N$ and
$1\le a<d':=N/d$ we have the distribution relation
\begin{equation}\label{equ:prodrel}
    \sum_{0\le j<d} L_N(1|a+jd')=  L_N(1|ad).
\end{equation}
It follows from the main result of Bass \cite{Bass} corrected by
Ennola \cite{En} that all the linear relations between $L_N(1|j)$
are consequences of \eqref{equ:symrel} and \eqref{equ:prodrel}.
Hence the seeded relations have the following forms in
words: for all $w\in \fA^0$
\begin{equation}\label{equ:seed}
\left\{\aligned
(N-2)Z(y_{1,j}*w-y_{1,-j}*w)=&(N-2j)(Z(y_{1,1}*w-y_{1,-1}*w),\\
\sum_{0\le j<d} Z(y_{1,a+jd'}*w)=&  Z(y_{1,ad}*w).
\endaligned\right.
\end{equation}

\section{Regularized distribution relations} \label{sec:Rac}
Multiple polylogs satisfy
the following distribution formula (cf.~\cite[Prop.~2.25]{Rac}):
\begin{equation} \label{equ:dist}
Li_{s_1,\dots,s_n}(x_1,\dots,x_n) = d^{s_1+\dots+s_n-n}
\sum_{y_j^d=x_j,1\le j\le n} Li_{s_1,\dots,s_n}(y_1,\dots,y_n),
\end{equation}
for all positive integer $d$. When $s_1=1$
we need to exclude the case of $x_1=1$. We call these \emph{finite
distribution relations} (FDT). Racinet
further considers the regularized version of these relations, which
we now recall briefly.

Fix an embedding  $\mmu_N\hookrightarrow \C$ and denote by $\gG$
its image. Define two sets of words
$$\bfX:=\bfX_\gG=\{x_\gs: \gs\in \gG\cup\{0\}\}, \quad\text{and}\quad
 \bfY:=\bfY_\gG=\{y_{n,\gs}=x_0^{n-1}x_\gs: n\in \N, \gs\in \gG\}.$$
Then one may consider the coproduct $\gD$ of $\Q\langle\bfX\rangle$ defined
by $\gD x_\gs=1\otimes x_\gs+x_\gs\otimes 1$ for all $\gs\in \gG\cup\{0\}$.
For every path $\gam\in \P^1(\C)-(\{0,\infty\}\cup\gG)$
Racinet defines the group-like element $\calI_\gam\in \CXX$ by
$$\calI_\gam:=\sum_{p\in \N,\gs_1,\dots,\gs_p\in  \gG\cup\{0\}}
 \calI_\gam(\gs_1,\dots,\gs_p)x_{\gs_1}\cdots x_{\gs_p},$$
where $\calI_\gam(\gs_1,\dots,\gs_p)$ is the iterated integral
$\int_\gam \om(\gs_1) \cdots \om(\gs_p)$ with
$$\om(\gs)(t)=\left\{
              \begin{array}{ll}
                \gs\,dt/(1-\gs t), & \hbox{if $\gs\ne 0$;} \\
                dt/t, & \hbox{if $\gs=0$.}
              \end{array}
            \right.
$$
This $I_\gam$ is essentially the same element denoted by $\dch$ in \cite{DG}.
Note that $\QY$ is the sub-algebra of $\QX$ generated by words not ending
with $x_0$. We let $\pi_\bfY:\QX\ra \QY$ be the projection. As $x_0$ is primitive
one knows that $(\QY,\gD)$ has a graded co-algebra structure.

Let $\QX_\cv$ be the sub-algebra of $\QX$ not beginning with $x_1$ and not
ending with $x_0$. Let $\pi_\cv:\QX\ra \QX_\cv$ be the projection. Passing
to the limit one get:
\begin{prop} \emph{(\cite[Prop.2.11]{Rac})} The series
$\calI_\cv:=\lim_{a\to 0^+,b\to 1^-}\pi_\cv(\calI_{[a,b]})$
is group-like in $(\CXX_\cv,\gD)$.
\end{prop}
Let $\calI$ be the unique group-like element in $(\CXX,\gD)$ whose coefficients
of $x_0$ and $x_1$ are 0 such that $\pi_\cv(\calI)=\calI_\cv$.
In order to do the numerical computation we need to find out explicitly the
coefficients for $\calI$. Put
$$\calI=\sum_{p\in \N, ,\gs_1,\dots,\gs_p\in  \gG\cup\{0\}}
C(\gs_1,\dots,\gs_p)x_{\gs_1}\cdots x_{\gs_p}.$$
\begin{prop} \label{prop:CoeffI}
Let $p,$ $m$ and $n$ be three non-negative integers. If
$p>0$ then we assume $\gs_1\ne 1$ and $\gs_p\ne 0$. Set
$(\gs_1,\dots,\gs_p,\{0\}^n) =(\gs_1,\dots,\gs_q)$. Then we have
 \begin{align} \notag
     \ C(\{1\}^m,\gs_1,& \dots,\gs_p,\{0\}^n)\\
     =& \left\{
          \begin{array}{ll}
            0, {\displaystyle \phantom{\frac1m}}  & \hbox{if $mn=p=0$;}  \\
            Z\big(\pi_\bfY(x_{\gs_1}\cdots x_{\gs_p})\big) , {\displaystyle \phantom{\frac1m}}  & \hbox{if $m=n=0$;}  \\
          {\displaystyle -\frac1m \sum_{i=1}^q C(\{1\}^{m-1},\gs_1,\dots,\gs_i,1,\gs_{i+1},\dots,\gs_q),} & \hbox{if $m>0$;} \\
                      {\displaystyle -\frac1n \sum_{i=1}^p C(\gs_1,\dots,\gs_{i-1},0,\gs_i,\dots,\gs_p,\{0\}^{n-1}),} & \hbox{if $m=0,n>0$.}
          \end{array}
        \right. \label{equ:CoeffI}
 \end{align}
Here $Z$ is defined by \eqref{equ:Z}.
\end{prop}
\begin{rem} This proposition provides the recursive relations we
may use to compute all the coefficients of $\calI$.
\end{rem}
\begin{proof} Since $\calI$ is group-like we have
\begin{equation}\label{equ:Igam}
\gD \calI=\calI\otimes\calI.
\end{equation}
The first case follows from this immediately since $C(0)=C(1)=0$.
The second case is essentially the definition \eqref{equ:Z} of $Z$.
If $m>0$ then we can compare the coefficient of
$x_1\otimes x_1^{m-1}x_{\gs_1}\cdots x_{\gs_q}$ of the two sides
of \eqref{equ:Igam} and find the relation \eqref{equ:CoeffI}.
Finally, if $m=0$ and $n>0$ then we may similarly consider the coefficient
of $x_{\gs_1}\cdots x_{\gs_p}x_0^{n-1}\otimes x_0$ in  \eqref{equ:Igam}.
This finishes the proof of the proposition.
\end{proof}

For any divisor $d$ of $N$ let $\gG^d=\{\gs^d:\gs\in \gG\}$,
$i_d:\gG^d\hookrightarrow \gG$ the embedding, and $p^d:\gG\twoheadrightarrow \gG^d$
the $d$th power map. They induce two algebra homomorphisms:
\begin{alignat*}{6}
p^d_\ast: \Q\langle\bfX_\gG\rangle & \lra\Q\langle\bfX_{\gG^d}\rangle & {} & &\hskip1cm
i_d^\ast: \Q\langle\bfX_\gG\rangle & \lra\Q\langle\bfX_{\gG^d}\rangle \\
x_\gs &\lms \begin{cases} d x_0, \  &\text{ if }\gs=0,\\
x_{\gs^d}, \  &\text{ if }\gs\in \gG.\end{cases}
 &\qquad\overset{\displaystyle\text{ and }}{\phantom{\sum}} & &
x_\gs &\lms \begin{cases} x_0, \  &\text{ if }\gs=0,\\
x_\gs, \  &\text{ if }\gs\in \gG^d,\\
0, \  &\text{ otherwise}.\end{cases}
\end{alignat*}

It is easy to see that both $i_d^\ast$ and $p^d_\ast$ are $\gD$-coalgebra morphisms
such that $i_d^\ast(\calI)$ and $p^d_\ast(\calI)$ have the same image under the map
$\pi_\cv$. By the standard Lie-algebra mechanism one has
\begin{prop} \label{prop:RDist} \emph{(\cite[Prop.2.26]{Rac})}
For every divisor $d$ of $N$
\begin{equation}\label{equ:RDist}
   p^d_\ast(\calI)=\exp\left(\sum_{\gs^d=1,\gs\ne 1} Li_1(\gs) x_1\right) i_d^\ast(\calI).
\end{equation}
\end{prop}
Combined with Proposition \ref{prop:CoeffI} the above result provides
the so-called \emph{regularized distribution relations} (RDT)
which of course include all the FDT of MPVs
given by \eqref{equ:dist}.

Computation suggests the the following conjecture concerning a
special class of distribution relations.
\begin{conj}\label{conj:distr}
Let $d$ be a positive integer. Then all the distribution relations in
\eqref{equ:dist}, where $x_j=1$ for all $j$, are consequences of
RDS of MPVs of level $d$.
\end{conj}
We are able to confirm this conjecture in the special case that
$w=2$, $n=1$, and $d$ is a prime.
\begin{thm} Write $L(i,j)=L_p(1,1|i,j)$ and $D(i)=L_p(2|i)$. Define
for $1\le i,j<p$:
\begin{align*}
\ & \text{FDT}:=
   D(0)-p\sum_{j=0}^{p-1}D(j)=0, \qquad
\text{RDS$(i)$}:=
 D(i)+L(i,0)-L(i,-i)=0,\\
\ &\text{FDS$(i,j)$}:=
 D(i+j)+L(i,j)+L(j,i)-L(i,j-i)-L(j,i-j)=0.
\end{align*}
Then
\begin{equation}\label{equ:distNull}
{\rm FDT}=\sum_{1\le i<p} {\rm FDS}(i,i)
+2\sum_{1\le j<i<p} {\rm FDS}(i,j)+2\sum_{i=1}^{p-1} {\rm RDS}(i,i).
\end{equation}
\end{thm}
\begin{proof} By changing the order of summation we see that
\begin{align*}
2\sum_{1\le j<i<p} D(i+j)=& \sum_{i=2}^{p-1}\sum_{j=1}^{i-1} D(i+j) +
\sum_{j=1}^{p-2}\sum_{i=j+1}^{p-1} D(i+j)\\
=&\sum_{i=2}^{p-2}\sum_{i\ne j=1}^{p-1} D(i+j)
+ \sum_{j=1}^{p-2} D(j-1) +\sum_{i=2}^{p-1} D(i+1)\\
=&(p-3)\sum_{j=0}^{p-1} D(j)-\sum_{i=2}^{p-2} D(i)-\sum_{i=1}^{p-1} D(2i)
+ \sum_{j=1}^{p-2} D(j) +\sum_{j=2}^{p-1} D(j)+2D(0)\\
=&(p-1)D(0)+(p-3)\sum_{j=1}^{p-1} D(j)
\end{align*}
since $\sum_{j=0}^{p-1} D(i+j)=\sum_{j=1}^{p-1} D(j)$ for all $i$ and
$\sum_{i=1}^{p-1} D(2i)=\sum_{i=1}^{p-1} D(i)$. This implies that the
dilogarithms on the right hand side of \eqref{equ:distNull} exactly add up to FDT.
Thus we only need to
show that all the double logarithms on the right hand side of \eqref{equ:distNull}
cancel.

First we note that $L(i,0)$ in ${\rm FDS}(i,i)$ and ${\rm RDS}(i,i)$ cancel. Now let
us consider the lattice points $(i,j)$ of $\Z^2$ corresponding to $L(i,j)$.
The points $(i,j)$ corresponding to $L(i,j)$ with positive signs fill in exactly
the area inside the square $[1,p-1]\times [1,p-1]$ (boundary inclusive):
$L(i,i)$ in ${\rm FDS}(i,i)$
provides the diagonal $y=x$, $\sum_{1\le j<i<p} L(i,j)$ (resp. $\sum_{1\le j<i<p} L(j,i)$)
form the lower right (resp. upper left) triangular region.

For the negative terms of the double logs, $L(i,-i)$ in ${\rm RDS}(i)$ provides the
diagonal $x+y=p$, $\sum_{1\le j<i<p} L(i,j-i)=\sum_{i=2}^{p-1}\sum_{j=p+1-i}^{p-1}L(i,j)$
form the upper right triangular region. Similarly, by  changing the order of summation
$\sum_{1\le j<i<p} L(j,i-j)=\sum_{i=1}^{p-2}\sum_{j=i+1}^{p-1}L(i,j-i)
=\sum_{i=2}^{p-1}\sum_{j=1}^{p-1-i}L(i,j)$ fills the lower left region.
\end{proof}

Further, numerical evidence up to level $N=49$ supports the following
\begin{conj}\label{conj:regdistr}
In weight two, all RDT are consequences of RDS and FDT.
\end{conj}

\section{Lifted relations from lower weights}
Note that when $N=3$ there are no seeded relations nor (regularized)
distribution relations.
When we deal with MZVs and alternating Euler sums we expect that
all the linear relations come from RDS.
Are these enough when $N=3$? Surprisingly, the answer is no.

The first counterexample is in weight four, i.e., $(w,N)=(4,3)$. Easy computation
shows that there are 144 MPVs in this case among which there are 239
nontrivial RDS which include 191 FDS.
Using these relations we get 127 independent linear
relations among the 144 MPVs. But the upper bound of $d(4,3)$ by
\cite[5.25]{DG} is 16, so there must be at least one more linearly
independent relation. Where else can we find it?
It is easy to verify that
all the seven RDT (including four FDT) can be derived from RDS.
However, we know that a product of two weight two MPVs
is of weight four. So on each of the five RDS (including two FDS)
in $\MPV(2,3)$ we can multiply any one of
the nine MPVs of $(w,N)=(2,3)$ to get a relation in $\MPV(4,3)$. For
instance, we have a FDS
$$ Z(y_{1,1}\ast y_{1,1}-y_{1,1}\sha y_{1,1})
    =L_3(2|2)+2L_3(1,1|1,1)-L_3(1,1|1,0)=0.$$
Multiplying by $L_3(1,1|1,1)=Z(y_{1,1}y_{1,2})$ we have a new relation
not derivable from RDS in $\MPV(4,3)$:
\begin{align*}
 & Z\big(y_{1,1}y_{1,2} \sha[y_{2,0}+2y_{1,1}y_{1,2}-2y_{1,1}y_{1,0}] \big)\\
=&L_3(1,1,2|1,1,0)+2L_3(1,2,1|1,1,0)+2L_3(2,1,1|1,1,0)+L_3(2,1,1|2,2,1)
+4L_3(\{1\}^4|1,1,2,1)\\
+&8L_3(\{1\}^4|1,0,1,0)-6L_3(\{1\}^4|1,0,0,1)
-4L_3(\{1\}^4|1,0,1,2)-2L_3(\{1\}^4|1,1,2,0)=0.
\end{align*}
Such relations coming from the lower weights are called
\emph{lifted relations (from lower weights).} In this way, when
$(w,N)=(4,3)$ we can produce 45 lifted RDS relations from weight two,
58 from weight three. We may also lift RDT and obtain
nine and six relations from weight two and three,
respectively. However, all the lifted relations together only
produce one new linearly independent relation, as expected.
Hence we find totally 128 linearly independent relations among the
144 MPVs of $(w,N)=(4,3)$. This implies that $d(4,3)\le 16$ which
is the same bound obtained by \cite[5.25]{DG} and is proved
to be exact under a variant of Grothendieck's period conjecture
by Deligne \cite{Del}.

For general levels $N$ we may lift not only RDS and RDT
but also the seeded relations. But a moment reflection tells us that
the lifted seeded relations are seeded so we don't need to
consider these after all.

\begin{defn}
We call a $\Q$-linear relation between MPVs \emph{standard}
if it can be produced by combinations of the following four families
of relations: regularized double shuffle relations (RDS),
regularized distribution relations (RDT), seeded
relations, and lifted relations from the above. Otherwise, it is
called a \emph{non-standard} relation.
\end{defn}

There are no seeded relations if $N=3$. In this case we believe
that all the linear relations among MPVs come from RDS and the lifted
relations (see Conjecture~\ref{conj:mainconj}). Moreover, computation
in small weight cases supports the following
\begin{conj} \label{conj:level3}
Suppose $N=3$ or $4$. Every MPV of level $N$ is a linear combination
of MPVs of the form $L(\{1\}^w|t_1,\dots,t_w)$ with $t_j\in
\{1,2\}$.  Consequently, the $\Q$-dimension of the MPVs of weight $w$
and level $N$ is given by $d(w,N)=2^w$ for all $w\ge 1$.
\end{conj}

\begin{rem}\label{rem:notenough}
Even adding all the lifted relations from lower weights does not
provide all the linear relations among MPVs. A quick look at the
Table~\ref{Ta:dbzeta} in \S\ref{sec:comp} tells us that if $(w,N)=(3,4)$
even though we know $d(3,4)\le 8$ and $d(4,4)\le 16$  by
\cite[5.25]{DG}, and the equality should hold by
Conjecture~\ref{conj:level3} or by a variant of Grothendieck's period conjecture
(see Remark~\ref{rem:Del}), we cannot produce enough relations
by using the standard ones. Instead, we can only show that
$d(3,4)\le 9$ and $d(4,4)\le 21$. More recently, by using octahedral
symmetry of $\P^1-(\{0,\infty\}\cup \mmu_4)$ we find (presumably all) the
non-standard relations in these two cases (see \cite{Zocta}).
\end{rem}

\begin{rem}\label{rem:Del}
Let $N=2,3,4$ or $8$. Assuming a variant of
Grothendieck's period conjecture, Deligne \cite{Del}
constructs explicitly a set of basis for $\MPV(w,N)$.
His results would imply that $d(w,2)$ is given by the Fibonacci numbers,
$d(w,3)=d(w,4)=2^w$, and $d(w,8)=3^w.$
\end{rem}

\section{Some conjectures of FDS and RDS}
Recall that if a map $Z_R:\fA^0\lra R$ satisfies the FDS and any
one of the equivalent conditions in Theorem~\ref{thm:RDS} then we
say that $Z_R$ has the \emph{regularized double shuffle} (RDS)
property. Let $R_{RDS}$ be the universal algebra (together with a
map $Z_{RDS}:\fA^0\lra R_{RDS}$) such that for every $\Q$-algebra
$R$ and a map $Z_R:\fA^0\lra R$ satisfying RDS there always exists
a map $\varphi_R$ to make the following diagram commutative:
\begin{equation*}
\text{$\diagramcompileto{diagwt2}
\fA^0 \drto_{Z_R}\rto^-{Z_{RDS}}& R_{RDS} \dto^-{\varphi_R}\\
\ & R
\enddiagram$}
\end{equation*}

When $N=3$ computation shows that the lifted relations contribute
non-trivially when the weight $w=5$: we can only get $d(5,3)\le 33$ instead
of the conjecturally correct dimension 32 without using lifted
relations. We may say that $Z_R$ has the
\emph{lifted regularized double shuffle} (LRDS)
property if it satisfies RDS and for all
$\om_1\in \fA^1$ and $\om_0,\om_0',\om_0''\in \fA^0$
$$Z_R(Z_R^{-1}\circ \rho_R\circ Z_R(\om_1)*\om_0-\om_1*\om_0)=Z_R(
 (\om_0*\om_0')*\om_0''-(\om_0\sha \om_0')*\om_0'')=0.$$
We can define $Z_{SR}$ and $R_{SR}$ corresponding to
the standard relations similar to  $Z_{RDS}$ and $R_{RDS}$ such that for
every $\Q$-algebra $R$ and a map $Z_R:\fA^0\lra R$ satisfying
the standard relations
there always exists a map $\varphi_R$ to make the following
diagram commutative:
\begin{equation}\label{equ:cm}
\text{$\diagramcompileto{diagwt2}
\fA^0 \drto_{Z_R}\rto^-{Z_{SR}}& R_{SR} \dto^-{\varphi_R}\\
\ & R
\enddiagram$}
\end{equation}

\begin{conj} \label{conj:mainconj}
Let $(R,Z_R)=(\R, Z)$ if $N=1,2$ and $(R,Z_R)=(\C, Z)$ if $N=p$ is a
prime $\ge 3$, where $Z$ is given by \eqref{equ:Z}.
If $N=1$ or $2$ then the map $\varphi_\R$ is injective, namely,
the algebra of MPVs is isomorphic to $R_{RDS}$. If $N=p$  is a
prime $\ge 3$ then the
map  $\varphi_\C$ is injective so the algebra of MPVs of level
$p$ is isomorphic to $R_{SR}$. Moreover, if $N=3$ then
 $Z_{LRDS}=Z_{SR}$ and $R_{LRDS}=R_{SR}$.
\end{conj}

 From Conjecture \ref{conj:mainconj} all the linear
relations among MPVs can be produced by RDS when $N=1$ or $2$,
and by the standard ones when $N=p$ is prime $\ge 3$. When $p\ge 5$
this is proved in \cite{Zocta} under the assumption of Grothendieck's
period conjecture.

Computation in many cases such as those listed in
Remark~\ref{rem:Ta:bN} and \ref{rem:wt2} show that MPVs must satisfy some other
relations besides the standard ones when $N$ has more than
two distinct prime factors, so a naive generalization of
Conjecture~\ref{conj:mainconj} to all levels does not exist at
present. However, when $N=4$ we find that octahedral symmetry of
$\P^1-(\{0,\infty\}\cup \mmu_4)$ may provide all the non-standard
relations (see \cite{Zocta}). But since we only have numerical evidence
in weight 3 and 4 it may be a little premature to form a conjecture
for level four at present.

\section{The structure of MPVs and some examples}\label{sec:comp}
In this section we concentrate on RDS between MPVs of
small weights. Most of the computations in this section are
carried out by MAPLE. We have checked the consistency of these
relations with many known ones and  verified our results numerically
using GiNac \cite{GiNac} and EZ-face \cite{EZface}.

By considering all the admissible words we see easily that the
number of distinct MPVs of weight $w\ge 2$ and level $N$ is
$N^2(N+1)^{w-2}$ and there are at most $N(N+1)^{w-2}$ RDS (but not FDS).
If $w\ge 4$ then the number of FDS is given by
$$(N-1)N^2(N+1)^{w-3}+\Big(\Big[\frac
w2\Big]-1\Big)N^4(N+1)^{w-4}=\Big(N^2\Big[\frac
w2\Big]-1\Big)N^2(N+1)^{w-4}.$$
If $w=2$ (resp.~$w=3$) then the number of FDS is $(N-1)^2$
(resp.~$N^2(N-1)$).

\subsection{Weight one.} \label{sec:wt1}
 From \S\ref{sec:seed} we know that all relations in weight one
follow from \eqref{equ:symrel} and \eqref{equ:prodrel}, and no RDS exists.
The relations in weight one are crucial for higher level cases
because they provide the seeded relations considered in \S\ref{sec:seed}.
Moreover, easy computation by
\eqref{equ:symrel} and \eqref{equ:prodrel} shows that there is
a hidden integral structure, namely, in each level there exists a
$\Q$-basis consisting of MPVs such that every other MPV is a $\Z$-linear
combination of the basis elements. This fact is proved by Conrad
\cite[Theorem~4.6]{Conrad}. Similar results should hold for higher
weight cases and we hope to return to this in a future publication
\cite{Zint}.

\subsection{Weight two.}
There are $N^2$ MPVs of weight 2 and level $N$:
 $$L_N(1,1|i,j),\quad L_N(2|j),   \qquad 1\le i\le N-1, 0\le j\le N-1.$$
For $1\le i,j<N$ the FDS $Z^*(y_{1,i}*y_{1,j})=Z^\sha(y_{1,i}\sha
y_{1,j})$ yields
  \begin{equation}\label{equ:wt2FDS}
L_N(2|i+j)+L_N(1,1|i,j)+L_N(1,1|j,i)=L_N(1,1|i,j-i)+L_N(1,1|j,i-j).
\end{equation}
Now from RDS $\rho(Z^*(y_{1,0}*y_{1,i}))= Z^\sha(y_{1,0}\sha
y_{1,i})$ we get for $1\le i<N$
\begin{equation}\label{equ:wt2RDS}
L_N(1,1|i,0)+L_N(2|i)=L_N(1,1|i,-i).
\end{equation}
The FDT in \eqref{equ:dist} yields: for every divisor $d$ of $N$,
and $1\le a,b< d':=N/d$
\begin{align}\label{equ:wt2dist1}
   L_N(2|ad)=&d\sum_{j=0}^{d-1}L_N(2|a+jd'), \\
   L_N(1,1|ad,bd)=&\sum_{j,k=0}^{d-1} L_N(1,1|a+jd',b+k d').\label{equ:wt2dist2}
\end{align}
To derive the RDT we can compare the coefficients of $x_1x_{\mu^{ad}}$
in \eqref{equ:RDist} and use Prop.~\ref{prop:CoeffI} to get:
for every divisor $d$ of $N$, and $1\le a< d'$
\begin{multline} \label{equ:wt2RDT}
    L_N(1|ad)\sum_{j=1}^{d-1}L_N(1|jd')
=\sum_{j=1}^{d-1}\sum_{k=0}^{d-1} L_N(1,1|jd',a+kd')\\
-\sum_{k=0}^{d-1} L_N(1,1|a+kd',-a-kd')-L_N(1,1|ad,-ad).
\end{multline}
By definition, the seeded relations are obtained from
\eqref{equ:symrel} and \eqref{equ:prodrel}. For example, if $N=p$ is
a prime then \eqref{equ:prodrel} is trivial and \eqref{equ:symrel}
is equivalent to: for all $1\le j<h:=(p-1)/2$
\begin{equation}\label{equ:symrelprime}
L_N(1|j)-L_N(1|-j)=(p-2j)(L_N(1|h)-L_N(1|h+1)).
\end{equation}
Thus multiplying by $L_N(1|i)$ ($1\le i<p$) and applying the shuffle
relation $L_N(1|a)L_N(1|b)=L_N(1,1|a,b-a)+L_N(1,1|b,a-b)$ we get:
\begin{multline}\label{equ:wt2N=pseed}
L_N(1,1|i,j-i)+L_N(1,1|j,i-j)-L_N(1,1|i,-j-i)-L_N(1,1|-j,i+j)\\
=(p-2j)\big(L_N(1,1|i,h-i)+L_N(1,1|h,i-h)
-L_N(1,1|i,-i-h)-L_N(1,1|-h,i+h) \big).
\end{multline}

Computation shows that the following conjecture should hold.
\begin{conj}
The RDT \eqref{equ:wt2RDT}
follows from the combination of the following relations:
the seeded relations, the RDS \eqref{equ:wt2FDS} and \eqref{equ:wt2RDS},
and the FDT \eqref{equ:wt2dist1} and \eqref{equ:wt2dist2}.
\end{conj}

\subsection{Weight three.}
Apparently there are $N^2(N+1)$ MPVs of weight 3 and level $N$:
for each choice $(i,j,k)$ with $1\le i\le N-1, 0\le j,k\le N-1$ we
have four MPVs of level $N$:
 $$L_N(1,1,1|i,j,k), \quad L_N(1,2|i,j), \quad L_N(2,1|j,k),  \quad L_N(3|k).$$
For $1\le i,j,k<N$ the FDS
$Z^*\big(y_{1,i}*(y_{1,j}y_{1,k})\big)=Z^\sha\big(y_{1,i}\sha(y_{1,j}y_{1,k})\big)$
yields
\begin{multline}\label{equ:wt3DS1}
L_N(\{1\}^3,i,j-i,k)+L_N(\{1\}^3,j,i-j,k+j-i)+L_N(\{1\}^3,j,k,i-k-j)\\
=L_N(2,1,i+j,k)+L_N(1,2,j,i+k)\hskip4cm\ \\
+L_N(\{1\}^3,i,j,k)+L_N(\{1\}^3,j,i,k)+L_N(\{1\}^3,j,k,i).
\end{multline}
For $1\le i,j<N$ the FDS $Z^*(y_{1,i}*y_{2,j})=Z^\sha(y_{1,i}\sha
y_{2,j})$ yields
\begin{equation}
  \label{equ:wt3DS2}
\aligned \ & L_N(3,i+j)+L_N(1,2,i,j)+L_N(2,1,j,i)\\
=&L_N(1,2,i,j-i)+L_N(2,1,i,j-i)+L_N(2,1,j,i-j).
\endaligned
\end{equation}
Moreover, there are three ways to produce RDS. Since $\gb(T)=T$
the first family of RDS come from
$Z^*\big(y_{1,0}*(y_{1,i}y_{1,i+j})\big)=Z^\sha\big(y_{1,0}\sha
      (y_{1,i}y_{1,i+j})\big)$ for $1\le i\le N-1, 0\le j\le N-1$:
\begin{align*}
\ &y_{1,0}*(y_{1,i}y_{1,i+j})
 =y_{1,0}y_{1,i}y_{1,i+j}+y_{1,i}
      \tau_i(y_{1,0}*y_{1,j})+y_{2,i}y_{1,i+j}\\
 =&y_{1,0}y_{1,i}y_{1,i+j}+y_{1,i}y_{1,i}y_{1,i+j}+y_{1,i}y_{1,i+j}y_{1,i+j}
        +y_{1,i}y_{2,i+j}+y_{2,i}y_{1,i+j}
\end{align*}
On the other hand,
 $$y_{1,0}\sha y_{1,i}y_{1,i+j}
      =y_{1,0}y_{1,i}y_{1,i+j}+y_{1,i}y_{1,0}y_{1,i+j}+y_{1,i}y_{1,i+j}y_{1,0}.$$
Hence
 \begin{multline}\label{equ:wt3RDS1}
 L_N(\{1\}^3|i,0,j)+L_N(\{1\}^3|i,j,0)+L_N(1,2|i,j)+L_N(2,1|i,j)\\
  =L_N(\{1\}^3|i,-i,i+j)+L_N(\{1\}^3|i,j,-i-j).$$
\end{multline}

The second family of RDS follow from $\gb(Z^*(y_{1,0}*y_{2,i}))=
Z^\sha(y_{1,0}\sha y_{2,i})$:
$$y_{1,0}y_{2,i}+y_{2,i}y_{1,i}+y_{3,i}=y_{1,0}y_{2,i}+y_{2,0}y_{1,i}+y_{2,i}y_{1,0}$$
which implies that
\begin{equation}\label{equ:wt3RDS2}
L_N(2,1,i,0)+L_N(3,i)=L_N(2,1,i,-i)+L_N(2,1,0,i).
\end{equation}

Now we consider the last family of RDS. By the definition of
stuffle product:
\begin{align*}
 y_{1,0}*y_{1,0}*y_{1,i}=&(2y_{1,0}^2+y_{2,0})*y_{1,i}\\
 =&2y_{1,0}(y_{1,0}*y_{1,i})+2y_{1,i}^3+2y_{2,i}y_{1,i}+y_{2,0}*y_{1,i} \\
 =&2y_{1,0}^2y_{1,i}+2y_{1,0}y_{1,i}^2+2y_{1,0}y_{2,i}+2y_{1,i}^3
          +2y_{2,i}y_{1,i}+y_{2,0}*y_{1,i}.
\end{align*}
Applying $\gb\circ Z^*$ and noticing that
$Z^\sha_{(2|0)}(T)=\zeta(2)$ we get
\begin{multline}
 \label{equ:Z*wt3}
(T^2+\zeta(2))Z^\sha_{(1|i)}(T)=
2Z^\sha_{(1,1,1|0,0,i)}(T)+2Z^\sha_{(1,1,1|0,i,i)}(T)+2Z^\sha_{(1,2|0,i)}(T)\\
+2Z^\sha_{(1,1,1|i,i,i)}(T)
          +2Z^\sha_{(2,1|i,i)}(T)+Z^\sha_{(2|0)}(T)Z^\sha_{(1|i)}(T).
\end{multline}
On the other hand by the definition of shuffle product
\begin{align*}
     \om_0\sha  \om_0\sha y_{1,i} =& 2\om_0^2 \sha \om_i
       =2 \om_0^2 \om_i+2 \om_0 \om_i \om_0+ 2\om_i \om_0^2\\
       =&2y_{1,0}^2y_{1,i}+2y_{1,0}y_{1,i}y_{1,0}+2y_{1,i}y_{1,0}^2
\end{align*}
Applying $Z^\sha$ we get
\begin{equation}
 \label{equ:Zshawt3}
T^2Z^\sha_{(1|i)}(T)=2Z^\sha_{(1,1,1|0,0,i)}(T)+2Z^\sha_{(1,1,1|0,i,0)}(T)
+2Z^\sha_{(1,1,1|i,0,0)}(T).
\end{equation}
We further have
\begin{align*}
 \ & Z^\sha(y_{1,0}y_{1,i}^2+y_{1,0}y_{2,i}-y_{1,0}y_{1,i}y_{1,0})\\
 =&Z^\sha{(1,1,1|0,i,i)}(T)+Z^\sha_{(1,2|0,i)}(T)-Z^\sha_{(1,1,1|0,i,0)}(T)\\
 =&2Z^\sha_{(1,1,1|i,0,0)}(T)-Z^\sha_{(2,1|i,0)}(T)-Z^\sha_{(2,1|0,i)}(T)
 -Z^\sha_{(1,1,1|i,0,i)}(T) -Z^\sha_{(1,1,1|i,i,0)}(T)
\end{align*}
where we have used the facts that
\begin{align*}
       Z^\sha_{(1,2|0,i)}(T)=&TZ^\sha_{(2|i)}(T)-Z^\sha_{(2,1|i,0)}(T)-Z^\sha_{(2,1|0,i)}(T)\\
       Z^\sha_{(1,1,1|0,i,i)}(T)=&TZ^\sha_{(1,1|i,i)}(T)-Z^\sha_{(1,1,1|i,0,i)}(T)
       -Z^\sha_{(1,1,1|i,i,0)}(T)\\
       Z^\sha_{(1,1,1|0,i,0)}(T)=&TZ^\sha_{(1,1|i,0)}-2Z^\sha_{(1,1,1|i,0,0)}(T)\\
      Z^\sha_{(1,1|i,0)} = &Z^\sha_{(2|i)}(T)+Z^\sha_{(1,1|i,i)}(T).
\end{align*}
Hence for $1\le i<N$ we have by subtracting \eqref{equ:Zshawt3}
from \eqref{equ:Z*wt3}
\begin{multline}\label{equ:wt3RDS3}
L_N(\{1\}^3|i,0,0)+L_N(2,1|i,0)+L_N(\{1\}^3|i,-i,0)=\\
L_N(2,1|i,-i)+L_N(2,1|0,i)+L_N(\{1\}^3|i,-i,i)+L_N(\{1\}^3|i,0,-i).
\end{multline}
Setting $j=0$ in \eqref{equ:wt3RDS1} and subtracting from
\eqref{equ:wt3RDS3} we get
 \begin{multline}
 \label{equ:wt3RDS4}
 L_N(\{1\}^3|i,-i,0)=
L_N(2,1|i,-i)+L_N(2,1|0,i)+L_N(\{1\}^3|i,0,0)+L_N(1,2|i,0).
\end{multline}

\subsection{Upper bound of $d(w,N)$ by Deligne and Goncharov}
By using the theory of motivic fundamental groups of
$\P^1-(\{0,\infty\}\cup \mmu_N)$ Deligne and Goncharov \cite[5.25]{DG} show that
$d(w,N)\le D(w,N)$ where $D(w,N)$ are defined by the formal power series
\begin{equation}\label{equ:DGbound}
 1+\sum_{w=1}^\infty D(w,N)t^w=
\left\{
    \begin{array}{ll}
      (1-t^2-t^3)^{-1}, & \hbox{if $N=1$;} \\
      (1-t-t^2)^{-1}, & \hbox{if $N=2$;} \\
     \big(1-\big(\frac{\varphi(N)}2+\nu(N)\big)t+
        \big(\nu(N)-1 \big)t^2 \big)^{-1}, & \hbox{if $N\ge 3$.} \\
    \end{array}
  \right.
\end{equation}
Here $\varphi$ is the Euler's totient function and $\nu(N)$
is the number of distinct prime factors of $N$.
Set $a=a(N):=\varphi(N)/2+\nu(N)$ and $b=b(N):=\nu(N)-1$. If $N>2$ then we have
$$ \sum_{j=1}^\infty D(w,N) t^n
 = at+(a^2-b)t^2+(a^3-2ab)t^3+(a^4-3a^2b+b^2)t^4+(3ab^2-4a^3b+a^5)t^5
 + \cdots$$
We will compare the bound obtained by standard relations to $D(w,N)$ in
the next two sections.

\section{Computational results in weight two}\label{sec:wt2comp}
In this section we combine the analysis in the previous sections
and the theory developed by Deligne and Goncharov \cite{DG} to present
a detailed computation in weight two and level $N\le 49$.

Let $\calG:=\iota(\Lie U_\om)$ be the motivic fundamental Lie algebra
(see \cite[(5.12.2)]{DG}) associated to the motivic fundamental group
of $\P^1-(\{0,\infty\}\cup \mmu_N)$. As pointed out in \S6.13 of op. cit. one may
safely replace $\calG(\mmu_N)^{(\ell)}$ by $\calG$ throughout \cite{G1}.
Then it follows from the proof of \cite[5.25]{DG} that
if a variant of Grothendieck's period conjecture (see 5.27(c) of op. cit.)
is true, which we assume in the following, then
\begin{equation}\label{equ:d(2,N)}
 d(2,N)=D(2,N)-\dim \ker(\gb_N),
\end{equation}
where
$\gb_N: \bigwedge^2 \calG_{-1,-1} \lra \calG_{-2,-2}$ is given by
Ihara's bracket $\gb_N(a \wedge b)=\{a,b\}$ defined by (5.13.6) of op. cit.
Here $\calG_{\bullet,\bullet}$
is the associated graded of the weight and depth
gradings of $\calG$ (see  \cite[\S 2.1]{G1}).
Let $k(N):=\dim \ker(\gb_N)$. Then
\begin{equation}\label{equ:delta1}
 \gd_1(N):=\dim \calG_{-1,-1}=\varphi(N)/2+\nu(N)-1
\end{equation}
by \cite[Thm.\ 2.1]{G1}. Thus
\begin{equation}\label{equ:iN}
 i(N):=\dim \im(\gb_N)=\gd_1(N)(\gd_1(N)-1)/2-k(N).
\end{equation}
Since $\dim \calG_{-2,-1}=\varphi(N)/2$ if $N>2$ and $0$ otherwise
the dimension of the degree 2 part of $\calG$ is
\begin{equation}\label{equ:delta2}
 \gd_2(N):=\dim \calG_{-2,-1}+\dim\calG_{-2,-2}=
\left\{
  \begin{array}{ll}
     i(N), & \hbox{if $N=1$ or 2;} \\
     \varphi(N)/2+i(N), & \hbox{if $N\ge 3$.}
  \end{array}
\right.
\end{equation}
Let $sr(N)$ be the upper bound of $\gd_2(N)$ obtained by
the standard relations. This can be computed by the method
described in \cite{Zocta}.  Let $SR(N)$ be the upper bound
of $d(2,N)$ similarly obtained by standard relations.
In Table \ref{Ta:dbzeta22} we use MAPLE to provide
the following data: $k(N)$, $sr(N)$, and  $SR(N)$. Then we can
calculate $\gd_1$, $\gd_2$ and $i(N)$ from \eqref{equ:delta1} to
\eqref{equ:delta2}. From \eqref{equ:d(2,N)}
we can check the consistency by verifying
$$sr(N)-\gd_2(N)=SR(N)-d(2,N)=SR(N)-D(2,N)+k(N)$$
which gives the number of linearly independent non-standard relations
(assuming Grothendieck's period conjecture). To save
space we use $D=D(2,N)$ and $d=d(2,N)$.

\begin{table}[h]
\begin{center}
\begin{tabular}{|c|c|c|c|c|c|c|c|c |c|c|c|c|c|c|c|c|c|c|c| }
 \hline
$N$      & 1 & 2 & 3 & 4 & 5 & 6 & 7 & 8 & 9 & 10 & 11 &12 &13 & 14 &15 & 16 & 17& 18  & 19\\ \hline
$\gd_1$  & 0 & 1 & 1 & 1 & 2 & 2 & 3 & 2 & 3 & 3  & 5  & 3 & 6 &  4 & 5 &  4 &  8 & 4  & 9  \\  \hline
$i$      & 0 & 0 & 0 & 0 & 0 & 1 & 1 & 1 & 3 & 3  & 5  & 3 & 8 &  6 &10 &  6 & 16 & 6  & 21 \\  \hline
$k$      & 0 & 0 & 0 & 0 & 1 & 0 & 2 & 0 & 0 & 0  & 5  & 0 & 7 &  0 & 0 &  0 & 12 & 0  & 15\\  \hline
$\gd_2$  & 0 & 1 & 1 & 1 & 2 & 2 & 4 & 3 & 6 & 5  & 10 & 5 & 14&  9 &14 & 10 & 24&  9  & 30\\  \hline
$sr$     & 0 & 1 & 1 & 1 & 2 & 2 & 4 & 4 & 6 & 6  & 10 & 8 & 14& 12 &16 & 16 & 24& 19  & 30\\  \hline
$D$      & 1 & 2 & 4 & 4 & 9 & 8 & 16& 9 & 16& 15 & 36 & 15& 49& 24 &35 & 25 & 81& 24  & 100\\  \hline
$SR$     & 1 & 2 & 4 & 4 & 8 & 8 & 14& 10& 16& 16 & 31 & 18& 42& 27 &37 & 31 & 69& 34  & 85\\  \hline
$d$      & 1 & 2 & 4 & 4 & 8 & 8 & 14& 9 & 16& 15 & 31 & 15& 42& 24 &35 & 25 & 69& 24  & 85\\  \hline
\end{tabular}
\begin{tabular}{ |c|c|c|c|c|c|c|c|c |c|c|c|c|c|c|c|c|} \hline
$N$      & 20 & 21 & 22 & 23 & 24 & 25 & 26 & 27 & 28 & 29 & 30 & 31 & 32 & 33 & 34 & 35 \\ \hline
$\gd_1$  & 5  & 7  & 6  & 11 &  5 & 10 & 7  & 9  & 7  & 14 & 6  & 15 &  8 & 11 & 9  & 13 \\  \hline
$i$      & 10 & 21 & 15 & 33 & 10 & 40 & 21 & 36 & 21 & 56 & 15 & 65 & 28 & 55 & 36 & 78 \\  \hline
$k$      & 0  &  0 &  0 & 22 &  0 &  5 &  0 &  0 &  0 & 35 &  0 & 40 &  0 & 0  & 0  & 0  \\  \hline
$\gd_2$  & 14 & 27 & 20 & 44 & 14 & 50 & 27 & 45 & 27 & 70 & 19 & 80 & 36 & 65 & 44 & 90 \\  \hline
$sr$     & 24 & 32 & 35 & 44 & 32 & 50 & 42 & 54 & 48 & 70 & 48 & 80 & 64 & 77 & 72 & 96 \\  \hline
$D$      & 35 & 63 & 48 &144 & 35 &121 & 63 &100 & 63 &225 & 47 &256 & 81 & 143& 99 & 195\\  \hline
$SR$     & 45 & 68 & 58 &122 & 53 &116 & 78 &109 & 84 &190 & 76 &216 &109 & 158& 127& 201\\  \hline
$d$      & 35 & 63 & 48 &122 & 35 &116 & 63 &100 & 63 &190 & 47 &216 & 81 & 143& 99 & 195\\  \hline
\end{tabular}
\begin{tabular}{ |c|c|c|c|c|c|c|c|c |c|c|c|c|c|c|} \hline
$N$      & 36 & 37 & 38 & 39 & 40 & 41 & 42 & 43 & 44 & 45 & 46 & 47 & 48 & 49 \\ \hline
$\gd_1$  & 7  & 18 & 10 & 13 &  9 & 20 & 8  & 21 & 11 & 13 & 12 & 23 &  9 & 21  \\  \hline
$i$      & 21 & 96 & 45 & 78 & 36 & 210& 28 &133 & 55 & 78 & 66 &161 & 36 & 175 \\  \hline
$k$      & 0  & 57 &  0 & 0  &  0 & 70 &  0 & 77 &  0 &  0 &  0 & 92 &  0 &  35  \\  \hline
$\gd_2$  & 27 & 114& 54 & 90 & 44 & 140& 34 &154 & 65 & 90 & 77 &184 & 44 & 196  \\  \hline
$sr$     & 72 & 114& 89 &112 & 96 & 140& 96 &154 &120 &144 &132 &184 &128 & 196  \\  \hline
$D$      & 63 & 361&120 &195 & 99 & 441& 79 &484 &143 &195 &168 &576 & 99 & 484 \\  \hline
$SR$     &108 & 304&156 &217 &151 & 371&141 &407 &198 &249 &223 &484 &183 & 449 \\  \hline
$d$      & 63 & 304&120 &195 & 99 & 371& 79 &407 &143 &195 &168 &484 & 99 & 449 \\  \hline
\end{tabular}
\caption{Upper bound of $d(2,N)$ obtained by standard relations and \cite[5.25]{DG}.}
\label{Ta:dbzeta22}
\end{center}
\end{table}

\begin{rem} \label{rem:wt2}
We now make the following observations in weight two case.

(a) If the level $N$ is a prime then the standard relations provide all
the $\Q$-linear relations under the assumption of a variant of
Grothendieck's period conjecture. This is proved in \cite[Thm.~1]{Zocta}.

(b) Notice that when $p\ge 11$ the vector space $\ker\beta_p$ contains a subspace
isomorphic to the space
of cusp forms of weight two on $X_1(p)$ which has dimension
$(p-5)(p-7)/24$ (see \cite[Lemma 2.3 \& Theorem 7.8]{G1}).
So it must contain another piece which
has dimension $(p-3)/2$ since
$\dim(\ker\beta_p)=(p^2-1)/24$ by \cite[(5)]{Zocta}. What is this missing piece?

(b) If $N$ is a $2$-power or a $3$-power then $D(2,N)$ should be sharp by
the conjecture mentioned in (a). See Remark~\ref{rem:Del}.

(c) If $N$ has at least two distinct
prime factors then $D(2,N)$ seems to be sharp,
though we don't have any theory to support it.

(d) Suppose Grothendieck's period conjecture is true. Then
by \cite[5.27]{DG}, (b) and (c) is equivalent to saying that the kernel
of $\gb_N$ is trivial if $N$ is a 2-power
or a 3-power, or has at least two distinct prime factors.
We believe this condition on $N$ for $\gb_N$ to be trivial
is necessary, too.

(e) If the level $N$ is a $p$-power for some prime
$p\ge 5$ then $\gb_N$ is unlikely to be injective (the
prime square case is proved in \cite[Prop.~5.3]{Zocta}).
We conjecture that non-standard relation doesn't exist
(i.e., $SR(N)$ is sharp), though we only have verified
the first two prime squares, $N=25$ and $N=49$.
\end{rem}

\begin{rem} \label{rem:Ta:bN}
In the three cases
$(w,N)=(2,8), (2,10)$ and $(2,12)$ we see that $d(w,N)>D(w,N)$.
By numerical computation we conjecture that the bounds given by
$D(w,N)$ are sharp in these cases and the following relations
are the non-standard ones: let $L_N(-)=L_N(1,1|-)$ and $L_N^{(2)}(-)=L_N(2|-)$,
then
{\allowdisplaybreaks
\begin{align}
37L_8(1,1)=&34L_8^{(2)}(5)+112L_8(3,1)+11L_8(3,0)+37L_8^{(2)}(1)-2L_8(2,6)\notag \\
  & +3L_8(7,3)-111L_8(5,7)+38 L_8(7,7)-8L_8(5,5),\label{equ:conj2}\\
 7L_{10}(5,2)=&72L_{10}^{(2)}(1)+265L_{10}^{(2)}(7)-7L_{10}(2,5)-467L_{10}(4,2)+467L_{10}(8,6)\notag \\
             &+14L_{10}(5,6)+64L_{10}(9,8)-164L_{10}(9,4)+166L_{10}(7,9)\notag \\
         & -260L_{10}(8,1)-66L_{10}(3,9) -7L_{10}(6,9)+7L_{10}(6,5).\label{equ:conj3}\\
L_{12}(8,7)=&5L_{12}^{(2)}(5)+8L_{12}(8,10)-6L_{12}(10,11)-8L_{12}(9,11)+L_{12}(10,9)\notag \\
        \ & -15L_{12}(8,1)+5L_{12}(9,10)+5L_{12}(6,1) -L_{12}(1,1)\notag \\
        \ &+6L_{12}(8,11)-11L_{12}(6,11)+8L_{12}(8,3)-L_{12}(11,8),\\
60L_{12}(8,11)=&38L_{12}(8,7)+348L_{12}(10,11)+502L_{12}(9,11)-492L_{12}(10,9)\notag \\
        \ & +600L_{12}(8,1)-552L_{12}(9,10)-154L_{12}(11,10)+20L_{12}(6,1)\notag \\
        \ & +261L_{12}(6,11)-502L_{12}(8,3)+221L_{12}(11,8)-319L_{12}(8,10),\\
221L_{12}(1,1)=&1854L_{12}(8,10)+562L_{12}(8,7)-1018L_{12}(10,11)-2416L_{12}(9,11)\notag \\
\ & +319L_{12}(10,9)-4270L_{12}(8,1)+2293L_{12}(9,10)+956L_{12}(11,10)\notag \\
\ &+1110L_{12}(6,1)+2416L_{12}(8,11)-3305L_{12}(6,11)+2416L_{12}(8,3).
\end{align}}
\end{rem}

\begin{defn}
We call the level $N$ \emph{standard} if either (i)
$N=1, 2$ or $3$, or (ii) $N$ is a prime power $p^n$ ($p\ge 5$).
Otherwise $N$ is called \emph{non-standard}.
\end{defn}

When $N$ is a non-standard level we find that very often
there are non-standard relations among MPVs.
For examples, the five relations in Remark~\ref{rem:Ta:bN} are
discovered only through numerical computation.
But are the standard relations enough to
produce all the linear relations when $N$ \emph{is} standard?
In weight two, when $N$ is a prime the answer is affirmative if
one assumes a variant of Grothendieck's period conjecture \cite{Zocta}.
Computations above give strong support for it and,
in fact, is the primary motivation of it.

\section{Computational results in other weights}\label{sec:wt3comp}
In this last section we briefly discuss our results in weight 3,4 and 5.
Since the computational complexity increases exponentially with weight
we cannot do as many cases as we have done in weight two.

Combining the FDS \eqref{equ:wt3DS1}, \eqref{equ:wt3DS2}, RDS
\eqref{equ:wt3RDS1}-\eqref{equ:wt3RDS4}, and the seeded relations
\eqref{equ:seed} we have verified the following facts by MAPLE:
$d(3,1)= 1, d(3,2)\le 3, d(3,3)\le 8$....
We have done similar computation in other small weight and low
level cases and listed the results in Table~\ref{Ta:dbzeta}.

We list some values of $D(w,N)$ in Table \ref{Ta:dbzeta} to compare with
the bound $SR(w,N)$ obtained by standard relations.
\begin{table}[h]
\begin{center}
\begin{tabular}{  |c|c|c|c|c|c|c|c |c|c|c|c|c|c| }
 \hline
$N $       & 1 & 2 & 3 & 4 & 5 & 6 & 7 & 8 & 9 & 10 & 11 & 12 & 13 \\ \hline
$SR(3,N)$  & 1 & 3 & 8 & 9 & 22& 23& 50& 38& 67& 70 & 157& 94 & 246 \\  \hline
$D(3,N)$   & 1 & 3 & 8 & 8 & 27& 21& 64& 27& 64& 56 & 216& 56 & 343 \\ \hline
$SR(4,N)$  & 1 & 5 & 16& 21& 61& 69&  &   &   & & &  &   \\  \hline
$D(4,N)$   & 1 & 5 & 16& 16& 81& 55& 256& 81 & 256& 209 & 1296 &209 & 2401  \\ \hline
$SR(5,N)$  & 2 & 8 & 32&   &   &   &   &   &   & &  & &   \\  \hline
$D(5,N)$   & 2 & 8 & 32& 32 & 243 & 144 & 1024& 243 & 1024& 780 & 7776 &780 & 16807  \\ \hline
\end{tabular}
\caption{Upper bound of $d(w,N)$ obtained by standard relations and \cite[5.25]{DG}.}
\label{Ta:dbzeta}
\end{center}
\end{table}
\begin{rem}
Note that $d(3,4)=D(3,4)+1$. By numerical computation we find the following non-standard relation:
\begin{align}
5L_4(1,2|2,3)=&46L_4(1,1,1|1,0,0)-7L_4(1,1,1|2,2,1)-13L_4(1,1,1|1,1,1) +13L_4(1,2|3,1) \notag \\
-L_4(1,1,1|&3,2,0)+25L_4(1,1,1|3,0,0)-8L_4(1,1,1|1,1,2) +18L_4(2,1|3,0), \label{equ:conj1}
\end{align}
Recently, we prove that by using the octahedral symmetry of
$\P^1-(\{0,\infty\}\cup \mmu_4)$ one can deduce equation \eqref{equ:conj1} (see \cite{Zocta}).
\end{rem}
 From the available data in Table \ref{Ta:dbzeta} we can
formulate the following conjecture.
\begin{conj} Let $N=p$ be a prime $\ge 5$. Then
$$d(3,p)\le \frac {p^3+4p^2+5p+14}{12}.$$
Moreover, equality hold if standard relations produce all the linear
relations.
\end{conj}
We obtained this conjecture under the belief that the upper
bound of $d(3,p)$ produced by the standard relations should
be a polynomial of $p$ of degree 3. Then we find the coefficients
by the bounds of $d(3,p)$ for $p=5,7,11,13$ in Table \ref{Ta:dbzeta}.

When $w>2$ it's not too hard to improve the bound of $d(w,p)$ given
in \cite[5.25]{DG} by the same idea as used in the proof
of \cite[5.24]{DG} (for example, decrease the bound by $(p^2-1)/24$).
But they are often not the best. We conclude our paper with the
following conjecture.
\begin{conj}\label{conj:stand}
If $N$ is a standard level then the standard relations always provide
the sharp bounds of $d(w,N)$, namely, all linear relations can be derived
from the standard ones. If $N$ is a non-standard level then
the bound in \emph{\cite[Cor.\ 5.25]{DG}} is sharp and the non-standard relations
exist in $\MPV(w,N)$ for all $w\ge 3$ (and in $\MPV(2,N)$ if $N\ge 10$).
\end{conj}

\end{document}